\documentclass{gtart_h}

\def\ifplaintex{\expandafter\ifx\csname documentclass\endcsname\relax}

\ifplaintex 
\else
\fi

\def\gt{{\mathsurround=0pt\it $\cal G\mskip-2mu$eometry \&\ 
$\cal T\!\!$opology}}        

\def\gtm{{\mathsurround=0pt\it $\cal G\mskip-2mu$eometry \&\ 
$\cal T\!\!$opology $\cal M\mskip-1mu$onographs}}    

\def\gtp{{\mathsurround=0pt\it $\cal G\mskip-2mu$eometry \&\ 
$\cal T\!\!$opology $\cal P\!$ublications}}  

\def\recd{{\small Received:\qua\receiveddate\ifx\reviseddate\relax
\else\qquad Revised:\qua\reviseddate\fi\par}} 

\def\volumenumber#1{\def\thevolumenumber{#1}}
\def\volumeyear#1{\def\thevolumeyear{#1}}
\def\volumename#1{\def\thevolumename{#1}}
\def\papernumber#1{\def\thepapernumber{#1}}
\def\pagenumbers#1#2{\def\startpage{#1}\def\finishpage{#2}}
\def\published#1{\def\publishdate{#1}}
\def\received#1{\def\receiveddate{#1}}
\def\revised#1{\def\reviseddate{#1}}
\def\accepted#1{\def\accepteddate{#1}}

\def\coverauthors#1{\def\thecoverauthors{#1}}
\def\asciiauthors#1{\def\theasciiauthors{#1}}
\def\asciiaddress#1{\def\theasciiaddress{#1}}
\def\asciiemail#1{\def\theasciiemail{#1}}

\def\coverauthors#1{\def\thecoverauthors{#1}}
\long\def\asciiabstract#1{\long\def\theasciiabstract{#1}}

\let\\\par
\let\thevolumenumber\relax\let\thepapernumber\relax
\let\thevolumeyear\relax\let\startpage\relax
\let\finishpage\relax\let\publishdate\relax\let\receiveddate\relax
\let\reviseddate\relax\let\accepteddate\relax\let\theasciititle\relax
\let\theasciiauthors\relax\let\theasciiaddress\relax
\let\theasciiabstract\relax
\let\thecoverauthors\relax
\let\thecoverauthors\relax\let\theerratum\relax\let\theasciiemail\relax
\let\theshortauthors\relax\let\theshorttitle\relax

\def\startpage{1}\def\finishpage{15}\def\thepapernumber{77}

\volumenumber{2}
\volumename{Proceedings of the Kirbyfest}
\volumeyear{1999}

\long\def\maketitlep{   

\count0=\startpage

\gtm\nl        
{\small Volume \thevolumenumber: \thevolumename\nl 
\ifx\theerratum\relax\else Erratum \erratumnumber\nl\fi
Pages \startpage--\finishpage\nl}

\vglue 0.1truein   

{\parskip=0pt\leftskip 0pt plus 1fil\def\\{\par\smallskip}{\ifplaintex\large
\else\Large\fi\bf\thetitle}\par\medskip}   
\vglue 0.05truein 

{\parskip=0pt\leftskip 0pt plus 1fil\def\\{\par}{\sc\theauthors}
\par\medskip}%
 
\vglue 0.03truein 

{\small\leftskip 25pt\rightskip 25pt{\bf Abstract}\stdspace\theabstract

{\bf AMS Classification}\stdspace\theprimaryclass
\ifx\thesecondaryclass\relax\else; \thesecondaryclass\fi\par
{\bf Keywords}\stdspace \thekeywords\par}\vglue 7pt

}   

\font\phead=cmsl9 scaled 950
\font=cmsl9 scaled 1050
\font\pnum=cmbx10 scaled 913
\font\lnum=cmbx10 
\font\pfoot=cmsl9 scaled 950
\font=cmsl9 scaled 1050
\ifplaintex
\headline{\vbox to 0pt{\vskip -4.5mm\line{\small\phead\ifnum
\count0=\startpage ISSN 1464-8997 (on line)
1464-8989 (printed) \hfill {\pnum\folio}\else\ifodd\count0\def\\{ }%
\ifx\theshorttitle\relax\thetitle\else\theshorttitle\fi\hfill{\pnum\folio}
\else\def\\{ and }{\pnum\folio}\hfill\ifx\theshortauthors\relax\theauthors
\else\theshortauthors\fi\fi\fi}\vss}}
\footline{\vbox to 0pt{\vglue 0mm\line{\small\pfoot\ifnum\count0=\startpage
Published \publishdate:\qua\copyright\ \gtp\hfill\else
\gtm, Volume \thevolumenumber\ (\thevolumeyear)\hfill\fi}\vss
}}
\else
\makeatletter
\def\@oddhead{{\small\lhead\ifnum\count0=\startpage ISSN 1464-8997 (on line)
1464-8989 (printed) \hfill {\lnum\number\count0}\else\ifodd\count0
\def\\{ }\ifx\theshorttitle\relax \thetitle \else\theshorttitle\fi\hfill
{\lnum\number\count0}\else\def\\{ and }{\lnum\number\count0}
\hfill\ifx\theshortauthors\relax 
\theauthors\else\theshortauthors\fi\fi\fi}}\def\@evenhead{@oddhead}
\def\@oddfoot{\small\lfoot\ifnum\count0=\startpage Published \publishdate:\qua\copyright\ \gtp\hfill\else
\gtm, Volume \thevolumenumber\ (\thevolumeyear)\hfill\fi}
\def\@evenfoot{@oddfoot}
\makeatother
\fi

\let\maketitlepage\maketitlep
\let\makeshorttitle\maketitlepage
\let\maketitle\maketitlepage

\newwrite\gtoutfile
\long\gdef\makeheadfile{  
{\def\\{, }\def\s{ }
\immediate\openout\gtoutfile head.xxx
\immediate\write\gtoutfile{Proxy-for: \ifx\theasciiauthors\relax
\theauthors\else\theasciiauthors\fi\s<\ifx\theasciiemail\relax\theemail\else\theasciiemail\fi>}
\immediate\write\gtoutfile{\noexpand\\}
\immediate\write\gtoutfile{Authors: \ifx\theasciiauthors\relax
\theauthors\else\theasciiauthors\fi}
{\def\\{ }\immediate\write\gtoutfile{Title: \ifx\theasciititle\relax
\thetitle\else\theasciititle\fi}}
\immediate\write\gtoutfile{Subj-class: GT or SG, GR etc}
\immediate\write\gtoutfile{MSC-class: \theprimaryclass\ifx\thesecondaryclass\relax\else, \thesecondaryclass\fi}
\immediate\write\gtoutfile{Journal-ref: Geom. Topol. Monogr. \thevolumenumber\s
(\thevolumeyear) \startpage-\finishpage}
\immediate\write\gtoutfile{Comments: Published by Geometry and Topology Monographs at}
\immediate\write\gtoutfile{\s\s\s  http://www.maths.warwick.ac.uk/gt/GTMon\thevolumenumber/paper\thepapernumber.abs.html}
\immediate\write\gtoutfile{\noexpand\\}
\immediate\write\gtoutfile{}
\ifx\theasciiabstract\relax
\immediate\write\gtoutfile{\theabstract}\else
\immediate\write\gtoutfile{\theasciiabstract}\fi
\immediate\write\gtoutfile{}
\immediate\write\gtoutfile{\noexpand\\}
\immediate\write\gtoutfile{}
\immediate\closeout\gtoutfile}}  

\def\maketitlepage{\maketitlep\makeheadfile}
\let\makeshorttitle\maketitlepage
\let\maketitle\maketitlepage

\volumenumber{7}
\volumename{Proceedings of the Casson Fest}
\volumeyear{2004}
\papernumber{7}
\pagenumbers{181}{203}
\received{30 December 2003}
\revised{28 April 2004}
\accepted{21 March 2004}
\published{18 September 2004}

\usepackage{amsmath,amssymb,amscd}

\newcommand\InjMod[1]{{\mathcal T}^+_{#1}}

\newcommand{\HF}{HF}

\newtheorem{theorem}{Theorem}[section]
\newtheorem{prop}[theorem]{Proposition}
\newtheorem{cor}[theorem]{Corollary}

\newtheorem{lemma}[theorem]{Lemma}

\newtheorem{defn}[theorem]{Definition}

\newcommand{\Q}{\mathbb{Q}}

\newcommand{\C}{\mathbb{C}}

\newcommand{\Z}{\mathbb{Z}}

\newcommand{\OneHalf}{\frac{1}{2}}

\newcommand{\Zmod}[1]{\Z/{#1}\Z}

\newcommand{\Image}{\operatorname{Im}}

\newcommand{\cm}{\cdot}

\newcommand{\ModSWfour}{\mathcal{M}}
\newcommand{\ModFlow}{\ModSWfour}

\newcommand{\SpinC}{{\mathrm{Spin}}^c}

\newcommand\abuts\Rightarrow
\newcommand\Sym{\mathrm{Sym}}

\newcommand\HFpRed{\HFp_{\red}}
\newcommand\HFpRedEv{\HFp_{\red, \ev}}
\newcommand\HFpRedOdd{\HFp_{\red,\odd}}
\newcommand\uHFpRedEv{\uHFp_{\red, \ev}}
\newcommand\uHFpRedOdd{\uHFp_{\red,\odd}}

\newcommand{\ev}{\mathrm{ev}}
\newcommand{\odd}{\mathrm{odd}}

\newcommand\relspinc{\underline{\spinc}}

\newcommand\x{\mathbf x}

\newcommand\z{\mathbf z}

\newcommand\y{\mathbf y}

\newcommand\ModSphere{\ModFlow\left({\mathbb S}\longrightarrow 
\Sym^{g-1}(\Sigma_{1})\times \Sym^2(\Sigma_{2})\right)}
\newcommand\ModSpheres\ModSphere
\newcommand\CF{CF}

\newcommand\CFa{\widehat{CF}}
\newcommand\CFp{\CFb}
\newcommand\CFm{\CF^-}

\newcommand\HFpEv{\HFp_{\mathrm{ev}}}
\newcommand\HFpOdd{\HFp_{\mathrm{odd}}}

\newcommand{\red}{\mathrm{red}}

\newcommand\HFp{\HFb}

\newcommand\HFm{\HF^-}
\newcommand\CFinf{CF^\infty}
\newcommand\HFinf{HF^\infty}
\newcommand\CFb{CF^+}
\newcommand\HFa{\widehat{HF}}
\newcommand\HFb{HF^+}

\newcommand\Mas{\mu}
\newcommand\UnparModSp{\widehat \ModSp}
\newcommand\UnparModFlow\UnparModSp
\newcommand\Mod\ModSp

\newcommand{\spinc}{\mathfrak s}

\newcommand\Real{\mathrm Re}

\newcommand\ModMaps{\mathcal M}
\newcommand\ModSp\ModMaps

\newcommand\Ta{{\mathbb T}_{\alpha}}
\newcommand\Tb{{\mathbb T}_{\beta}}

\newcommand\alphas{\mbox{\boldmath$\alpha$}}

\newcommand\betas{\mbox{\boldmath$\beta$}}

\newcommand\uHF{\underline{\HF}}
\newcommand\uHFp{\underline{\HF}^+}

\newcommand\uHFpRed{\uHFpred}
\newcommand\uHFpred{\underline{\HF}^+_{\red}}

\newcommand\uHFinf{\uHF^\infty}

\newcommand\Dual{\mathcal D}
\newcommand\Duality\Dual

\newcommand\spincrel\relspinc

\newcommand\CFK{CFK}
\newcommand\HFK{HFK}

\newcommand\CFKinf{\CFK^{\infty}}

\newcommand\HFKa{\widehat\HFK}

\title{On Heegaard Floer homology and\\Seifert fibered surgeries}

\author{Peter Ozsv\'ath\\Zolt{\'a}n Szab{\'o}}
\asciiauthors{Peter Ozsvath\\Zoltan Szabo}
\coverauthors{Peter Ozsv\noexpand\'ath\\Zolt{\noexpand\'a}n Szab{\noexpand\'o}}
\address{Department of Mathematics, Columbia University, NY 10025, USA\\
  Institute for Advanced Study, Princeton, NJ 08540, USA}
\secondaddress{Department of Mathematics, Princeton University, NJ 08544,
USA}
\asciiaddress{Department of Mathematics, Columbia University, 
NY 10025, USA\\Institute for Advanced Study, Princeton, 
NJ 08540, USA\\and\\Department of Mathematics, Princeton University, NJ
08544, USA}

\gtemail{\mailto{petero@math.columbia.edu}, \mailto{szabo@math.princeton.edu}}
\asciiemail{petero@math.columbia.edu, szabo@math.princeton.edu}

\begin{document}

\begin{abstract}  
  We explore certain restrictions on knots in the three--sphere which
  admit non-trivial Seifert fibered surgeries. These restrictions stem
  from the Heegaard Floer homology for Seifert fibered spaces, and
  hence they have consequences for both the Alexander polynomial of
  such knots, and also their knot Floer homology.  In particular, we
  show that certain polynomials are never the Alexander polynomials of
  knots which admit homology three--sphere Seifert fibered surgeries.
  The knot Floer homology restrictions, on the other hand, apply also
  in cases where the Alexander polynomial gives no information, such
  as the Kinoshita--Terasaka knots.
\end{abstract}

\asciiabstract{%
We explore certain restrictions on knots in the three-sphere which
  admit non-trivial Seifert fibered surgeries. These restrictions stem
  from the Heegaard Floer homology for Seifert fibered spaces, and
  hence they have consequences for both the Alexander polynomial of
  such knots, and also their knot Floer homology.  In particular, we
  show that certain polynomials are never the Alexander polynomials of
  knots which admit homology three-sphere Seifert fibered surgeries.
  The knot Floer homology restrictions, on the other hand, apply also
  in cases where the Alexander polynomial gives no information, such
  as the Kinoshita-Terasaka knots.}

\maketitle
\section{Introduction}

It is an interesting open question to characterize knots $K$ in the
three--sphere with the property that some (non-trivial) surgery on
$S^3$ along $K$ is a Seifert fibered space. This question has received
considerable attention recently, cf \cite{Dean,Gordon,EudaveMunoz}.
The aim of the present article is to present some
obstructions for a knot $K$ to admitting such surgeries.  These
obstructions in turn come from the Heegaard Floer homology
of \cite{HolDisk} and the related knot invariant defined
in \cite{Knots} and \cite{RasmussenThesis}.

\subsection{Seifert fibered surgeries and the Alexander polynomial}

For simplicity, throughout most of this paper we consider the case of
homology sphere surgeries, although other surgeries are accessible to
Floer homology as well (compare also
\cite[Subsection 8.1]{KMOS}). Specifically,
given a rational number $r$ and a knot $K\subset S^3$, let $S^3_r(K)$
denote the three--manifold obtained as Dehn surgery on $S^3$ along $K$.
Writing $r=p/q$ as a fraction in its lowest terms,
$H_1(S^3_{p/q}(K);\Z)\cong \Zmod{p}$. In particular, if $S^3_r(K)$ is
an integer homology three--sphere, then $r$ has the form $1/q$ for a
non-zero integer $q$.

We give a certain obstruction to Seifert fibered surgeries. The
strongest form of this result should be stated in terms of the
Heegaard Floer homology of the zero--surgery of $S^3$ along $K$,
$\HFp(S^3_0(K))$ cf Theorem \ref{thm:HFpObstruction} and
Proposition \ref{prop:NegativeCase} below.
However, for the purposes of this introduction, we prefer to state a
weaker form in terms of its Euler characteristic, which can be
expressed purely in terms of the Alexander polynomial. (This result is
proved in Section \ref{sec:Surgeries}.)

\begin{theorem}
\label{intro:AlexanderObstruction}
Let $K\subset S^3$ be a knot in the three--sphere. Write its
symmetrized Alexander polynomial as
$$\Delta_K(T)=a_0 + \sum_{i>0}a_i(T^i+T^{-i}),$$
and let
\begin{equation}
\label{eq:DefTorsion}
t_i(K)=\sum_{j=1}^\infty j a_{|i|+j}
\end{equation}
denote the torsion coefficients of the knot.
Then, if there is
an integer $q\neq 0$ for which $S^3_{1/q}(K)$ is Seifert fibered,
all the non-zero integers $t_i(K)$ have the same sign.
\end{theorem}

For example, for 34 of the 54 non-alternating knots with fewer than
eleven crossings, there are both positive and negative torsion
coefficients. It follows from Theorem \ref{intro:AlexanderObstruction}
that these knots admit no Seifert fibered $1/q$--surgeries. (Of course,
there are other known techniques for ruling out Seifert fibered surgeries for
sufficiently small knots, cf \cite{SnapPea}.)

This result rules out many alternating knots, as well. In particular,
there is an interesting family which we describe below.  But first,
recall that the sign of the $t_i(K)$ is governed by the following
result, proved in \cite{AltKnots}:
\begin{theorem}[Corollary 1.6
  of \cite{AltKnots}]
\label{thm:AltKnotsEstimate}
Let $K$ be an alternating knot, and let 
$$\delta(m,i)=\max\left(0,\left\lceil\frac{|m|-2|i|}{4}\right\rceil\right)$$
(note that this is the $i^{th}$ torsion coefficient of the $(2,2m+1)$
torus knot).  Then,
$$(-1)^{i+\frac{\sigma}{2}}(t_i(K)-\delta(\sigma,i))\leq 0,$$
where here $\sigma$ denotes the signature of $K$.
\end{theorem}
We now have the following result:
\begin{cor}
\label{cor:AltCase}
If $K$ is an alternating knot with genus $g$ and whose signature $\sigma$ 
satisfies
\begin{eqnarray*}
g+\frac{\sigma}{2}\equiv 1 \pmod{2}&{\text{and}}&
\sigma\neq 0,
\end{eqnarray*}
then there is no  non-zero
integer $q\neq 0$ for which $S^3_{1/q}(K)$ is Seifert fibered.
\end{cor}
\begin{proof}
It is a classical result that the genus $g$ of an alternating knot
agrees with the degree of its Alexander polynomial
(cf \cite{Crowell,MurasugiAlt}). In particular, $t_{g-1}(K)\neq 0$.

Recall that in general $|\sigma|\leq 2g$, while our
hypothesis on the parity of $\sigma/2$ forces this inequality to be strict,
and hence
$\delta(\sigma,g-1)=0$.  Theorem \ref{thm:AltKnotsEstimate}, together
with the hypothesis on the parity of $\sigma/2$ again, gives 
\begin{equation}
\label{eq:OneTorsion}
t_{g-1}(K)<0.
\end{equation}
On the other hand, letting $i=\frac{|\sigma|}{2}-1$, we have that
$\delta(\sigma,i)=1$, and hence another application of Theorem
\ref{thm:AltKnotsEstimate} gives
\begin{equation}
\label{eq:AnotherTorsion}
1\leq t_{i}.
\end{equation}
Together, inequalities \eqref{eq:OneTorsion} and \eqref{eq:AnotherTorsion},
combined with Theorem \ref{intro:AlexanderObstruction}, give the claimed result.
\end{proof}

To elucidate the hypotheses of the above corollary, note that the
figure--eight knot satisfies the parity hypothesis, but $\sigma=0$.
However, $+1$--surgery on this knot
gives the Brieskorn sphere $\Sigma(2,3,7)$.  Moreover,
for each integer $n\geq 1$, the $(2,2n+1)$ torus knot has non-zero
signature ($2n$), but $\sigma/2=g$. Of course, $+1$ surgery on any of
these knots is Seifert fibered.

The sign of the torsion coefficients as in
Theorem \ref{intro:AlexanderObstruction} is dictated by the
orientation on the Seifert fibered space. Specifically, note that for
any integral homology Seifert fibered space $Y$ different from $S^3$,
Casson's invariant $\lambda(Y)$ is non-zero (cf \cite{FSseifert,
NeumannWahl,FukuharaMatsumotoSakamoto}). Moreover,
Casson's surgery formula shows that
$$\lambda(S^3_{1/q}(K))=q\cm \sum_{i} t_i(K).$$
Putting these two together, we see that the signs
of the non-vanishing $t_i$ are determined by Casson's invariant of
$Y$. We say that $Y$ a Seifert fibered integral homology three--sphere
has a {\em positive Seifert orientation} if its Casson's invariant is
positive (ie has the same sign as that of the Brieskorn sphere
$-\Sigma(2,3,5)$).  In Section \ref{sec:Surgeries} we give a more
direct description of this sign.

\subsection{Knot Floer homology}

We have also a Seifert fibered surgery obstruction, which can be
stated in terms of the ``knot Floer homology''. This invariant is a
refinement of Heegaard Floer homology of the three--manifold $Y$, in
the presence of an oriented knot $K$ in $Y$, \cite{Knots,
RasmussenThesis}. We will say little about this invariant
beyond its formal properties. If $K$ is a null-homologous knot in a
three--manifold $Y$, then there is an induced filtration of the
Heegaard Floer complex, and the homology of the associated graded
object is this knot Floer homology.

For knots $K$ in the three--sphere, the invariants
$\HFKa$ take the form of graded Abelian groups indexed by integers
$i$.  According to a ``skein exact sequence'' which these
groups satisfy, it follows that if $\Delta_K(T)$ denotes the
symmetrized Alexander polynomial of $K$, then
\begin{equation}
\label{eq:EulerKnot}
\Delta_K(T)=\sum_{i\in\Z}\chi(\HFKa(K,i))\cm
T^i\in\Z[T,T^{-1}]
\end{equation}
(see \cite[Proposition 4.2]{Knots}). 

Moreover, according to
\cite[Theorem 1.2]{HolDiskGenus}, the knot
Floer homology determines the genus of the knot:
\begin{equation}
\label{eq:DetectsGenus}
g(K)=\max\{i\big|\HFKa(K,i)\neq 0\}.
\end{equation}
Knot Floer homology is difficult to calculate in general. However,
there are some families of knots for which the answer is known. For
example, if $K$ is an alternating knot, then $\HFKa(K,i)\cong
\Z^{|a_i|}$ is supported in dimension $i-\frac{\sigma}{2}$, where
$\sigma$ denotes the signature of $K$
(cf \cite[Theorem 1.4]{AltKnots}, see
also \cite{Rasmussen} for the two--bridge case). Here, of course, $a_i$
is the $T^i$ coefficient of the Alexander polynomial of $K$.

\subsection{Seifert fibered surgeries and knot Floer homology}

Theorem \ref{thm:HFpObstruction}, together with some additional
properties of the knot Floer homology which we review in
Section \ref{sec:KnotInvariants}, give strong restrictions on
$\HFKa(K)$ for knots which admit Seifert fibered surgeries. This
restriction leads to the following:

\begin{theorem}
\label{intro:HFKObstruction}
Let $K\subset S^3$ be a knot with genus $g$.  
If there is an integer $q>0$ so that $S^3_{1/q}(K)$
is a positively oriented Seifert fibered space, then $\HFKa(S^3,K,g)$
is trivial in odd degrees (and non-trivial in even degrees).
If there is an integer $q>0$ so that $S^3_{1/q}(K)$
is a negatively oriented Seifert fibered space and $g>1$,
then $\HFKa(S^3,K,g)$ is trivial in even degrees (and
non-trivial in odd degrees).
\end{theorem}

This has the following corollary (compare also \cite{KMOS}):

\begin{cor}
\label{cor:AlexHFKDefect}
If $\deg\Delta_K<g(K)$, then for integers $q\geq 0$,
$S^3_{1/q}(K)$ is never a positively oriented Seifert fibered space.
Indeed, if in addition,  $g(K)>1$ 
then no $1/q$ surgery along $K$ is Seifert fibered.
\end{cor}

\begin{proof}
  Note that $S^3_{1/q}(K)=-S^3_{-1/q}(r(K))$, where $r(K)$ denotes the
  reflection of $K$.  Thus, the corollary follows immediately from the Euler
  characteristic relation (equation \eqref{eq:EulerKnot}), 
  together with Theorem \ref{intro:HFKObstruction}.
\end{proof}

As an illustration, we consider the family of Kinoshita--Terasaka knots
$KT_{r,n}$ (see \cite{KinoshitaTerasaka}) with $|r|\geq 2$ and $n\neq
0$ (to avoid the unknot). These knots all have trivial Alexander
polynomial. However, Gabai has shown \cite{GabaiKT} that $KT_{r,n}$
has genus $|r|$ (note that one could deduce this alternatively from
the calculation of $\HFKa(K_{r,n},r)$ from \cite{calcKT}).  Thus,
Corollary \ref{cor:AlexHFKDefect} applies: none of these knots admits
surgeries which are integral homology
Seifert fibered spaces. See \cite{calcKT} for
other families of knots with these properties.

In this paper, we have dealt with restrictions on Seifert fibered surgeries
which are obtained with the use of Heegaard Floer homology.
Further results on Seifert fibered surgeries using monopole Floer
homology can be found in \cite[Subsection 8.1]{KMOS}.

\subsection{Surgeries giving $\Sigma(2,3,5)$ and $\Sigma(2,3,7)$}
\label{subsec:SpecialCases}

It is conjectured that if $K$ is a knot in $S^3$ which admits a
surgery which is $\Sigma(2,3,5)$, then $K$ is the trefoil
(cf remarks following \cite[Problem 3.6(D)]{Kirby}, see also
\cite{Zhang}), and if $K$ is a knot in $S^3$ which admits a surgery
which is $\Sigma(2,3,7)$, then $K$ is the figure--eight or the trefoil.
We use the surgery long exact sequence to provide some evidence for
this conjecture in the following form:

\begin{theorem}
\label{thm:PoincareSphere}
Let $K$ be a knot in $S^3$ with the property that 
$S^3_{r}(K)\cong \Sigma(2,3,5)$ for some $r\in \Q$, 
then $r=-1$, and the knot Floer homology of $K$ agrees with
that of the left-handed trefoil. In particular, the Seifert genus
of $K$ is one and its Alexander polynomial is $T^{-1}-1+T$.
\end{theorem}

\begin{theorem}
\label{thm:Bries2}
Let $K$ be a knot in $S^3$ with the property that $S^3_r(K)\cong
\Sigma(2,3,7)$ (as an oriented manifold)
then there are two cases. Either $r= -1$, in which case $K$ has genus
one and its Alexander polynomial agrees with that of the trefoil.  In
the case where $r=+1$, $K$ has genus one and its Alexander
polynomial is $-T^{-1}+3-T$.
\end{theorem}

\subsection{Acknowledgements}
The authors would like to thank Cameron Gordon for helpful discussions
on these questions.  The first author was partially supported by NSF
grant numbers DMS--0234311, DMS--0111298, and FRG--0244663, and the second
author was partially supported by NSF grant numbers DMS--0107792
and FRG--0244663, and a Packard Fellowship.

\section{A brief review of Heegaard Floer homology}
\label{sec:Background}

Very briefly, the Heegaard Floer homology is a homology theory
associated to a closed, oriented three--manifold $Y$.  For simplicity,
we start with the case where $Y$ is an integral homology three--sphere.
To define Heegaard Floer homology, we start with a suitably generic
Heegaard diagram $(\Sigma,\alphas,\betas,z)$, where here $\Sigma$ is
an oriented Riemann surface, and $\alphas=\{\alpha_1,...,\alpha_g\}$
and $\betas=\{\beta_1,...\beta_g\}$ are $g$--tuples of (pairwise
disjoint, embedded) attaching circles for the two handlebodies, and
$z$ is a point which does not lie on any of the attaching circles.  We
now consider the $g$--fold symmetric product $\Sym^g(\Sigma)$, ie
unordered $g$--tuples of points (counted with multiplicity) on
$\Sigma$, together with the pair of embedded tori
\begin{eqnarray*}
\Ta=\alpha_1\times ...\times \alpha_g
&{\text{and}}&\Tb=\beta_1\times...\times \beta_g.
\end{eqnarray*}

Loosely speaking, Heegaard Floer homology measures an obstruction to
pulling apart these two tori. More precisely, we fix a complex structure on 
$\Sigma$, which in turn induces a complex structure on its $g$--fold symmetric
product (with respect to which the two tori are totally real), and we consider
the chain complex $\CFinf(Y)$
generated by $[\x,i]\in (\Ta\cap \Tb)\times \Z$, and whose differential
counts holomorphic Whitney disks. That is, 
fix intersection points $\x, \y\in \Ta\cap\Tb $. A {\em Whitney disk} $u$
connecting $\x$ to $\y$ is a map 
$$u\co \{\z\in\C\big| |z|\leq 1\} \longrightarrow \Sym^g(\Sigma)$$
satisfying the boundary conditions
\begin{eqnarray*}
u\{\zeta \big| \Real(\zeta)\geq 0~{\text{and}}~|\zeta|=1\}\subset \Ta,
&&
u\{\zeta\big| \Real(\zeta)\leq 0~{\text{and}}~|\zeta|=1\}\subset \Tb,\\
u(-\sqrt{-1})=\x, &&
u(\sqrt{-1})=\y.
\end{eqnarray*}
The space of homotopy classes of Whitney disks connecting $\x$ to $\y$
is denoted $\pi_2(\x,\y)$. For a fixed Whitney disk $u$, let $n_z(u)$
denote the algebraic intersection number of $u$ with the
submanifold $\{z\}\times \Sym^{g-1}(\Sigma)\subset \Sym^g(\Sigma)$. 
Note that $n_z(u)$ depends
only on the homotopy class $\phi$ of $u$.

With these preliminaries in place, we can define a map
$$\partial[\x,i] = \sum_{\y\in\Ta\cap\Tb}
\sum_{\phi\in\pi_2(\x,\y)\big|\Mas(\phi)=1} c(\phi)\cm
[\y,i-n_z(\phi)],$$ where here $\Mas(\phi)$ denotes the expected
dimension of the moduli space of pseudo-holomorphic representatives for 
the homotopy class $\phi$, and 
$c(\phi)$ denotes an appropriate signed count of points
in this moduli space, modulo the natural
automorphism group. Here the term ``pseudo-holomorphic'' means a
sufficiently generic perturbation of the notion of holomorphic disks.
By adapting arguments from Lagrangian Floer homology (which in turn
rest on Gromov's compactness theorem \cite{Gromov}), one can show that
$\partial^2=0$. It is simplest to think of $c(\phi)$ as an element of
a field (rather than working over $\Z$), which we fix now to be the
field $\Q$ of rational numbers.

In fact, since $\{z\}\times \Sym^{g-1}(\Sigma)$ is a subvariety of
$\Sym^g(\Sigma)$ which is disjoint from $\Ta$ and $\Tb$, it follows
that if a given homotopy class $\phi$ has holomorphic representatives,
then $n_z(\phi)\geq 0$. In particular, it follows that the subset of
$\CFm(Y)\subset \CFinf(Y)$ generated by pairs $[\x,i]$ with $i<0$ is actually a
subcomplex.
Indeed, we have a short exact sequence of complexes
\begin{equation}
\label{eq:HFminfp}
0 \longrightarrow \CFm(Y) \longrightarrow \CFinf(Y)
  \stackrel{{\pi}}{\longrightarrow} \CFp(Y) \longrightarrow 0,
\end{equation}
where here $\CFp(Y)$ is defined to make the sequence exact; it is
generated by pairs $[\x,i]$ with $i\geq 0$.
All three complexes are endowed with an endomorphism, defined by
$U[\x,i]=[\x,i-1]$, which can be used to construct a fourth variant:
\begin{equation}
\label{eq:HFapp}
0 \longrightarrow \CFa(Y) \longrightarrow \CFp(Y)
  \stackrel{{U}}{\longrightarrow} \CFp(Y) \longrightarrow 0,
\end{equation}
where here $\CFa(Y)$ is, by definition, the kernel of the induced
endomorphism.  According to \cite{HolDisk}, the homology groups of the
complexes $\CFm(Y)$, $\CFinf(Y)$, $\CFp(Y)$, and $\CFa(Y)$, denoted
$\HFm(Y)$, $\HFinf(Y)$, $\HFp(Y)$, and $\HFa(Y)$, are topological
invariants of the three--manifold $Y$. In particular, they are
independent of the perturbations, complex structures, and Heegaard
diagrams which go into their definitions. In addition, the complexes
can be given a canonical $\Zmod{2}$--grading, characterized by the
properties that it is preserved by the maps in the short exact
sequences~\eqref{eq:HFminfp} and~\eqref{eq:HFapp}, and also that
$\chi(\HFa(Y))=1$, bearing in mind that we have assumed that
$H_1(Y;\Z)=0$. Sometimes, we write $\HFa_{\ev}(Y)$ (respectively
$\HFa_{\odd}(Y)$) for the subgroup of $\HFa(Y)$ generated by elements
with even (respectively odd) $\Zmod{2}$--grading (and the corresponding
notation for the other types of Heegaard Floer homology).

As a special case, consider $S^3$. For this manifold, $\HFa(S^3)\cong \Q$,
$\HFm(S^3)\cong \Q[U]$, $\HFinf(S^3) \newline \cong \Q[U,U^{-1}]$, and
$\HFp(S^3) \cong \Q[U,U^{-1}]/U\cm \Q[U]$ (thought of as a $\Q[U]$--module).
We denote this latter module by $\InjMod{}$. It is interesting to 
note that, according to
\cite[Theorem 10.1]{HolDiskTwo},
$$\HFinf(Y)\cong\Q[U,U^{-1}]$$ 
(supported in even parity)
for any integer homology three--sphere $Y$.

In the case of a general closed, oriented three--manifold the same
constructions work, except that when $b_1(Y)>0$, we use a restricted
class of Heegaard diagram, a technical point which we will not pursue
here (cf \cite[Section 5]{HolDisk}). But
there is some additional structure which is important to us here: the
chain complex splits into summands indexed by $\SpinC$ structures over
$Y$. For example, consider the 
case where $Y$ is obtained by $p$--framed surgery on a knot $K\subset S^3$,
where $p$ is an integer. 
In this case, there is an identification
$H_1(S^3_p(K);\Z)\cong
\Zmod{p}$. And indeed, we have a corresponding
identification $\SpinC(S^3_p(K))\cong \Zmod{p}$. 
In this case, there is a splitting of all three of the 
Heegaard Floer homologies indexed by integers $i\in\Zmod{p}$; we write
\begin{equation}
\label{eq:UniversalCoefficients}
\HFp(S^3_p(K))\cong \bigoplus_{i\in\Zmod{p}}\HFp(S^3_p(K),i).
\end{equation}
(The canonical $\Zmod{2}$--grading is slightly more complicated 
when $b_1(Y)>0$, 
but can still be defined,
cf \cite[Section 10.4]{HolDiskTwo},
see also equation \eqref{eq:HFinfTwist} below.)

In general, $\HFp(Y)$ is infinitely generated. However the quotient
$\HFpRed(Y)=\HFp(Y)/\HFinf(Y)$ is finitely generated. The link between
Heegaard Floer homology and the Alexander polynomial is provided by
the following result, which can be seen as an analogoue of the
Meng--Taubes theorem in Seiberg--Witten theory \cite{MengTaubes, TuraevSW}.

\begin{theorem}
\label{thm:Euler}
Let $K\subset S^3$ be a knot, and let $S^3_0(K)$ denote the
three--manifold obtained by performing $0$--framed surgery on $S^3$
along $K$. Then, $\HFp(S^3_0(K),i)$ is finitely generated when $i\neq
0$, with Euler characteristic given by
$$\chi(\HFp(S^3_0(K),i))=-t_i(K).$$ In the case where $i=0$,
$$\chi(\HFpRed(S^3_0(K),i))\geq -t_0(K).$$
\end{theorem}

The above theorem is stated and proved in
\cite[Theorem 5.2]{HolDiskTwo} in the case
where $i\neq 0$, while in the case where $i=0$, it is 
\cite[Theorem 10.17]{HolDiskTwo}.

A key calculational device we will use here is the following integer
surgeries long exact sequence, which is a restatement of
\cite[Theorem 9.19]{HolDiskTwo}:

\begin{theorem}
\label{thm:IntegerSurgeries}
Let $K\subset S^3$, and fix an integer $p>0$. Then, there is an 
long exact sequence
$$\ldots \longrightarrow \HFp(S^3)  \stackrel{{a}}{\longrightarrow}
  \bigoplus_{j\equiv i\!\!\!\!\pmod{p}}\HFp(S^3_0(K),j)
  \stackrel{{b}}{\longrightarrow} \HFp(S^3_p(K),[i]) \longrightarrow
\ldots,$$
where here $a$ and $c$ preserves the canonical $\Zmod{2}$ degree,
while $b$ reverses it.
\end{theorem}

There are many variants of this sequence. The case where $p=1$
generalizes to the case of arbitrary framed knots in an arbitrary
closed, oriented three--manifold (though the statements about the
$\Zmod{2}$ degree have to be suitably modified). (We say little about
the maps on Floer homology in the exact sequence, except that they are
induced from cobordisms. For more on this, see \cite{HolDiskFour}.)

Another variant uses fractional surgeries. The term associated to the
zero--surgery is slightly different in this case. Specifically, suppose
that $Y$ is a three--manifold with $H_1(Y;\Z)\cong \Z$, fix a
homomorphism $H^1(Y;\Z)\longrightarrow \Zmod{p}$. Then, there is a
notion of Floer homology with coefficients in $\Q[\Zmod{p}]$, written
$$\uHFp(Y;\Q[\Zmod{p}])\cong
\bigoplus_{i\in\Z}\uHFp(Y,i;\Q[\Zmod{p}]).$$ In fact, this theory is
obtained from a ``universal'' case of coefficients twisted by
$\Q[H^1(Y;\Z)]=\Q[\Z]$ which, of course, can be thought of as the ring
of Laurent polynomials with coefficients in $\Q$ over a formal
variable $T$. In this case, the universal coefficients theorem
applies, to show that $$\HFp_{i}(Y;\Q) \cong
\left(\uHFp_i(Y;\Q[\Z])\otimes \Q\right) \oplus 
\left(\Q\star \uHFp_{i+1}(Y;\Q[\Zmod{p}]))\right),$$
where here $i$ denotes the parity $i\in\Zmod{2}$, 
tensor products are taken over the base ring $\Q[\Z]$, where here
$\Q$ is viewed as a module over
$\Q[\Z]$ with trivial action by $\Z$, and
$A\star B$ denotes the torsion product over the principal
ideal domain $\Q[\Z]$.
More generally (ie generalizing from the case where $p=1$), we have
\begin{multline*}
\uHFp_{i}(Y;\Q[\Zmod{p}]) \cong \\
\left(\uHFp_i(Y;\Q[\Z])\otimes \Q[\Zmod{p}]\right) \oplus 
\left(\Q[\Zmod{p}]\star \uHFp_{i+1}(Y;\Q[\Zmod{p}]))\right),
\end{multline*}
where here $\Q[\Zmod{p}]$ is a module over $\Q[\Z]$ using our given
homomorphism $\Z\cong H^1(Y;\Z) \longrightarrow \Zmod{p}$. Note that
all of the groups $\HFa$, $\HFm$, $\HFinf$, and $\HFp$ have their
analogues with twisted coefficients, and they are related by exact
sequences analogous to equations \eqref{eq:HFminfp}
and \eqref{eq:HFapp}. In particular, we can once again form
the group 
$$\uHFpRed(Y;\Q[\Zmod{p}])=\uHFp(Y;\Q[\Zmod{p}])/\uHFinf(Y;\Q[\Zmod{p}]).$$ 
Again the canonical $\Zmod{2}$--grading gives a decomposition
$$\uHFpRed(Y;\Q[\Zmod{p}])=\uHFpRedEv(Y;\Q[\Zmod{p}])\oplus
\uHFpRedOdd(Y;\Q[\Zmod{p}]).$$
The Euler characteristic of $\uHFp$ with coefficients in $\Q[\Zmod{p}]$ is
related to the untwisted case by the relation 
\begin{equation}
\label{eq:TwistedChi}
        \chi\left(\uHFp(Y,i,\Q[\Zmod{p}])\right)= 
        p\cm \chi\left(\HFp(Y,i)\right)
\end{equation}
for all $i\neq 0$
(cf \cite[Lemma 11.1]{HolDiskTwo}).

We  use twisted coefficients 
for the following variant of the surgery long exact sequence,
which is a restatement of
\cite[Theorem 9.14]{HolDiskTwo}

\begin{theorem}
\label{thm:ExactFrac}
Let $K\subset S^3$ be a knot. For each integer $p>0$, there is a
homomorphism $H_1(S^3_0(K);\Z)\longrightarrow \Zmod{p}$ for which we
have a long exact sequence
$$\ldots  \longrightarrow \HFp(S^3)
  \stackrel{{a}}{\longrightarrow} \uHFp(S^3_0(K),\Q[\Zmod{p}])
  \stackrel{{b}}{\longrightarrow} 
  \HFp(S^3_{1/p}(K)) \longrightarrow \ldots,$$
where once again $a$ and $c$ preserve $\Zmod{2}$ degree and $b$ reverses it.
\end{theorem}

\subsection{Absolute gradings}

In the case where $Y$ is an integral homology three--sphere, the
$\Zmod{2}$--grading can be lifted to a $\Z$--grading on the Floer
homology of $\HFp(Y)$, cf
\cite[Section~7.1]{HolDiskFour}.  This provides
at once a numerical invariant for integer homology three--spheres, the
correction term denoted $d(Y)$, which is the minimal dimension of any
homogeneous element in $\HFp(Y)$ coming from $\HFinf(Y)$,
cf \cite{AbsGraded}.

This invariant restricts the intersection forms of smooth
four--manifolds which bound $Y$, according to the following result
(compare also the analogous gauge-theoretic results of
Fr{\o}yshov, cf \cite{FroyInst, FroySW}):

\begin{theorem} (Corollary 9.8 of \cite{AbsGraded})
\label{thm:NoIntForm}
If $Y$ is an integer homology three--sphere which bounds a smooth,
negative-definite four--manifold,
then $d(Y)\geq 0$.
\end{theorem}

In the case where $Y$ is not an integer homology three--sphere, the
absolute grading can be defined only for the summands of the Floer
homology which belong to $\SpinC$ structures over $Y$ whose first
Chern class is torsion. Moreover, the absolute grading is no longer a
$\Z$--grading in general, but rather, it gives a grading by rational
numbers \cite{AbsGraded}. We do not discuss this in great generality,
but rather content ourselves with the situation of three--manifolds $Y_0$ with
first homology isomorphic to $\Z$. 

In this case, let $\HFp(Y_0,0)\subset \HFp(Y_0)$ denote the summand
corresponding to the $\SpinC$ structure with vanishing first Chern
class.  In this case, the absolute grading takes its values in
rational numbers of the form $\OneHalf + \Z$. The relationship between
this absolute grading and the $\Zmod{2}$--grading referred to earlier
is given by the relation that elements whose absolute $\Q$--grading
lies in $-\OneHalf + 2\Z$ have even parity, while those whose absolute
$\Q$--grading lies in $\OneHalf + 2\Z$ have odd parity. 

For such three--manifolds, $\HFinf$ is determined by
\cite[Theorem 10.1]{HolDiskTwo} to have the
form
$$\HFinf(Y_0,0)\cong \Q[U,U^{-1}]\oplus\Q[U,U^{-1}],$$
where here the two summands have different parity. Correspondingly,
there are now two numerical invariants analogous to the earlier correction
term, $d_{\pm 1/2}(Y_0)$ -- the maximal $\Q$--grading of any element
in $\HFp(Y_0)$ contained in the image of $\HFinf(Y_0)$ whose parity is
given by $\pm \OneHalf + 2\Z$.

These numerical invariants can be used to give a more precise form of
Theorem \ref{thm:Euler} in the case where $i=0$. Specifically, we get
the following straightforward combination of
Theorem 10.17 with the discussion from
\cite[Section 4]{AbsGraded}:

\begin{equation}
\label{eq:RefinedEuler}
\chi(\HFpRed(S^3_0(K),0))-
\left(\frac{d_{-1/2}(S^3_0(K))-d_{1/2}(S^3_0(K))+1}{2}\right)
=-t_0(K)
\end{equation}

Note that it follows from the algebraic structure of $\HFinf$ (specifically,
an action of $H_1(Y_0;\Z)$ on $\HFinf(Y_0;\Z)$) that
\begin{equation}
\label{eq:DIneq}
d_{1/2}(Y_0)-1\leq d_{-1/2}(Y_0)
\end{equation}
(cf \cite[Proposition 4.10]{AbsGraded}).

\subsection{More on twisted coefficients}

Throughout this section, $Y_0$ will denote a three--manifold with
$H_1(Y_0;\Z)\cong \Z$, and indeed we implicitly fix such an
identification. This at once gives rise to an isomorphism
$\SpinC(Y_0)\cong \Z$ (via the map $\spinc \mapsto c_1(\spinc)/2$,
followed by the identification $H^2(Y_0;\Z) \cong \Z$ induced by
Poincar\'e duality), and we write $\HFp(Y_0,i)$ for the summand of
$\HFp(Y_0)$ in the $\SpinC$ structure corresponding to the integer
$i$. Similarly, there is an induced surjective homomorphism
$H^1(Y_0;\Z)\longrightarrow \Zmod{p}$. 

The following fact will be useful in the applications of
Section \ref{sec:KnotInvariants}.

\begin{prop}
\label{prop:TwistedGrows}
The rank of
$\uHFpRedEv(Y_0;\Q[\Zmod{p}])$,
and respectively that of
$\uHFpRedOdd(Y_0;\Q[\Zmod{p}])$, is a non-decreasing function
of $p$.
In fact, there is 
there is a short exact sequence 
\begin{equation}
\label{eq:TwistedGrows}
0 \longrightarrow  C \longrightarrow
  \uHFpRed(Y_0)\otimes_{\Q[T,T^{-1}]} \Q[\Zmod{p}]
  \longrightarrow \uHFpRed(Y_0;\Q[\Zmod{p}]) \longrightarrow 0,
\end{equation}
where $C$ is a $\Zmod{2}$--graded group of rank $(d_{-1/2}-d_{1/2}+1)/2$
supported entirely in odd parity.
\end{prop}

Before giving the proof, we prove another useful lemma
(which ensures that the $\Q$--rank of the 
vector space in the middle of Equation~\eqref{eq:TwistedGrows} grows
linearly with $p$). Note
$\uHFp(Y_0;\Q[T,T^{-1}])$ is a direct sum of cyclic
$\Q[T,T^{-1}]$--modules. As such, it has a canonical torsion submodule whose
quotient is a free $\Q[T,T^{-1}]$--module. 

\begin{lemma}
\label{lemma:TwistedTorsion}
The torsion submodule of $\uHFp(Y_0;\Q[T,T^{-1}])$ is in fact
annihilated by $1-T$. Indeed, this torsion submodule is identified with
the image of $\uHFinf(Y_0;\Q[T,T^{-1}])$ in
$\uHFp(Y_0;\Q[T,T^{-1}])$, and hence its quotient
$\uHFpRed(Y_0;\Q[T,T^{-1}])$ is a free $\Q[T,T^{-1}]$--module.
\end{lemma}

\begin{proof}
This can be seen, for example, from the twisted long exact sequence.
First, note that we can realize $Y_0$ as $0$--framed surgery along a
knot $K$ in an integer homology three--sphere $Y$. This gives rise
to an integral surgeries long exact sequence with twisted coefficients
\cite[Theorem 9.23]{HolDiskTwo},
generalizing Theorem \ref{thm:IntegerSurgeries}. In particular, when $p$
is sufficiently large, we get an exact sequence
\begin{multline*}
\ldots
 \stackrel{{\underline a}}{\longrightarrow} 
\uHFp(Y_0,i;\Q[\Zmod{p}])
 \stackrel{{\underline b}}{\longrightarrow} 
\HFp(Y_p,[i])\otimes_\Q \Q[\Zmod{p}] \\
 \stackrel{{\underline c}}{\longrightarrow} 
\HFp(Y)\otimes_\Q \Q[\Zmod{p}]
 \longrightarrow \ldots
\end{multline*}
Note that
\cite[Theorem 9.23]{HolDiskTwo} 
is stated for the zero--surgery in its ``universally twisted'' case,
i.e. with coefficients in $\Z[\Z]$, but this can be readily specialized
to $\Q[\Zmod{p}]$ as above; also, in general, the group corresponding
to the zero--surgery consists of a direct sum of groups 
$\uHFp(Y_0,j;\Q[\Zmod{p}])$ over all $j\equiv i\pmod{p}$, but 
here we apply this in the case where $p$ is sufficiently large that
each of these terms has only one non-zero summand.

The image of ${\underline a}$ is easily seen to coincide
with the image of $\uHFinf(Y_0;\Q[\Zmod{p}])$ inside
$\uHFp(Y_0;\Q[\Zmod{p}])$; while the image of ${\underline b}$
is a free $\Q[\Zmod{p}]$--module (since it is a submodule of
a free $\Q[\Zmod{p}]$--module).
\end{proof}

According to \cite[Theorem 10.12]{HolDiskTwo},
\begin{equation}
\label{eq:HFinfTwist}
\uHFinf(Y_0,0;\Q[\Z]) \cong \Q[U,U^{-1}],
\end{equation}
thought of as a $\Q[\Z]$--module with trivial action by $\Z$. (Indeed,
a result of this type holds for all three--manifolds and this is what
is used to define the canonical $\Zmod{2}$--grading in general.)  From
the universal coefficients theorem, it follows readily now that
$$\uHFinf(Y_0,0;\Q[\Zmod{p}]) \cong \Q[U,U^{-1}]\oplus \Q[U,U^{-1}]$$
for all $p>0$, where here the two summands have different parity.  It
might appear that this gives a new pair of 
correction terms  $d_{\pm 1/2}(Y_0;\Q[\Zmod{p}])$. However,
we have the following:

\begin{lemma}
\label{lemma:IndepTwist}
The correction term $d_{\pm 1/2}(Y;\Q[\Zmod{p}])$ is independent of $p$
(ie it agrees with the correction term $d_{\pm 1/2}(Y)$).
\end{lemma}

\begin{proof}
We have $\Q[T,T^{-1}]$--module maps
$$\Q \stackrel{{\eta}}{\longrightarrow}  \Q[\Zmod{p}]
  \stackrel{{\epsilon}}{\longrightarrow}  \Q$$
where $\eta(1) = \sum_{i=0}^{p-1} T^i$ and $\epsilon(T^i)=1$, so that
$\epsilon\circ \eta$ is multiplication by $p$ (ie an isomorphism).
In the notation, we use here an identification $\Q[\Zmod{p}]\cong
\Q[T]/T^p-1$.  These homomorphisms induce maps on the Floer homologies
with twisted coefficients; and indeed we get a diagram:
\begin{equation}
\label{eq:DIdent}
\begin{CD}
\HFinf(Y) @>{\eta^\infty}>{\cong}> \uHFinf(Y;\Q[\Zmod{p}]) @>{\epsilon^\infty}>{\cong}> \HFinf(Y) \\
@V{\pi}VV @V{\underline\pi}VV @V{\pi}VV \\
\HFp(Y) @>{\eta^+}>> \uHFp(Y;\Q[\Zmod{p}]) @>{\epsilon^+}>> \HFp(Y).
\end{CD}
\end{equation}
On the chain level, $\epsilon^+\circ \eta^+$ is multiplication by $p$, and hence
the induced maps on homology compose to give an isomorphism. Indeed, by
the universal coefficients theorem, both ${\eta^\infty}$ and ${\epsilon^\infty}$ induce isomorphisms. The statement about the correction terms now
follows from a diagram chase of Equation~\eqref{eq:DIdent}.
\end{proof}

Now, we prove Proposition \ref{prop:TwistedGrows}.

{\bf Proof of Proposition~\ref{prop:TwistedGrows}}\qua
The universal coefficients theorem gives
a diagram
$$
\begin{tiny}
\begin{CD}
0\longrightarrow\uHFinf(Y_0;\Q[\Z])\otimes \Q{[\Zmod{p}]}
@>>>\uHFinf(Y_0;\Q[\Zmod{p}])@>>>
\Q[\Zmod{p}]\star\uHFinf(Y_0;\Q[\Z])
\longrightarrow 0 \\
 @V{\alpha}VV @V{\beta}VV @VV{\gamma}V \\
0\longrightarrow\uHFp(Y_0;\Q[\Z])\otimes \Q{[\Zmod{p}]}
@>>>\uHFp(Y_0;\Q[\Zmod{p}])@>>>
\Q[\Zmod{p}]\star\uHFp(Y_0;\Q[\Z])
\longrightarrow 0 \\
 @VVV @VVV &  \\
 \uHFpRed(Y_0)\otimes\Q[\Zmod{p}] @>{f}>>
\uHFpRed(Y_0;\Q[\Zmod{p}]) \\
\end{CD}
\end{tiny}
$$
We claim first that $f$ is surjective. This follows from the fact that
$\gamma$ is surjective, a fact which we establish using the following
diagram
$$\begin{tiny}
\begin{CD}
0@>>>\Q[\Zmod{p}]\star\uHFinf(Y_0;\Q[\Z])
@>>>
\uHFinf(Y_0;\Q[\Z])@>{1-T^p\equiv 0}>> \uHFinf(Y_0;\Q[\Z]) \\
&& @V{\gamma}VV @VV{\underline\pi}V @VVV \\
0@>>>\Q[\Zmod{p}]\star\uHFp(Y_0;\Q[\Z]))
@>>>\uHFp(Y_0;\Q[\Z])
@>{1-T^p}>>\uHFp(Y_0;\Q[\Zmod{p}]),
\end{CD}
\end{tiny}$$
since according to Lemma \ref{lemma:TwistedTorsion}, the kernel of
$1-T$ (on the bottom row) is contained in the image of ${\underline
\pi}$.

Now, the kernels of $\alpha$, $\beta$, and $\gamma$ and the maps
between them are independent of $p$: they are determined by $d_{\pm
1/2}(Y_0;\Q[\Zmod{p}])$, which we saw in Lemma \ref{lemma:IndepTwist}
to be independent of $p$. In particular, the cokernel of the induced
map from $\ker \beta$ to $\ker \gamma$ (which in turn is identified with
the kernel of $f$) is a group supported entirely in odd parity, whose
rank is independent of $f$. Indeed, it is not hard to see
that this rank is given by $(d_{-1/2}(Y_0)-d_{1/2}(Y_0)+1)/2$. 
This provides us with the claimed
exact sequence from equation \eqref{eq:TwistedGrows}, establishing the
proposition.
\endproof

It is worth noting the following quick consequence of the above proposition:
\begin{multline}
\label{eq:RefinedEulerTwist}
\chi(\uHFpRed(Y_0,0;\Q[\Zmod{p}])) \\ -
\left(\frac{d_{-1/2}(Y_0;\Q[\Zmod{p}])-d_{1/2}(Y_0;\Q[\Zmod{p}])+1}{2}\right)
= -p\cm t_0(K).
\end{multline}

\section{Relationship with Heegaard Floer homology for the zero--surgery}
\label{sec:Surgeries}

Let $Y$ be a Seifert fibered space with $b_1(Y)=0$ or $1$. Such a
manifold can be realized as the boundary of a four--manifold
$W(\Gamma)$ obtained by plumbing two--spheres according to a weighted
tree $\Gamma$. Here, the weights are thought of as a map $m$ from the
set of vertices of $\Gamma$ to $\Z$.

\begin{defn}
\label{def:PositiveSeifert}
Suppose that $\Gamma$ is a weighted tree which has either
negative-definite or negative-semi-definite intersection form.  Then,
we say that the induced orientation 
$-\partial W(\Gamma)$ is a {\em positive Seifert orientation}.
\end{defn}

Note that if $Y$ is a Seifert fibered space with $b_1(Y)=0$ then at
least one of $+Y$ or $-Y$ has a positive Seifert orientation.
Moreover, either orientation on any lens space is a positive Seifert
orientation; similarly, either orientation on a Seifert fibered space
with $b_1(Y)=1$ is a positive Seifert orientation.  Finally, if $Y$ is
the quotient of a circle bundle $\pi\co N\longrightarrow \Sigma$
over a Riemann surface by a finite group of orientation-preserving
automorphisms $G$, and if $N$ is oriented as a circle bundle with
positive degree, then the induced orientation on $Y$ is a positive
Seifert orientation. An integral homology Seifert fibered space $Y$
has positive Seifert orientation if and only if its Casson invariant
is positive, cf \cite{FSseifert}.

Following \cite{SomePlumbs}, we work with the following convenient
generalization of the notion of positively oriented Seifert fibered
spaces:

\begin{defn}
  Let $\Gamma$ be a weighted graph which is a disjoint union of trees.
  The degree of a vertex $d(v)$ is the number of edges which contain
  $v$.  We say that $\Gamma$ is a {\em negative-definite (respectively
    negative-semi-definite) graph with at most one bad point} if the
  intersection form for $\Gamma$ is negative-definite (respectively
  negative-semi-definite), and there is at most one vertex
  $v\in\Gamma$ whose weight $m(v)$ is larger than $-d(v)$.
\end{defn}

For three--manifolds $Y=-\partial W(\Gamma))$, where $\Gamma$ is a
negative-definite graph with at most one bad point, $\HFp(Y)$ can be
explitly calculated in terms of the graph $\Gamma$,
see \cite{SomePlumbs}. Indeed, the part of that calculation which we
will use in the present paper can be summarized as follows:

\begin{theorem}
\label{thm:SeifertCalc}
Let $\Gamma$ be a negative-definite graph with at most one bad point.
Then, $\HFp(-\partial W(\Gamma))$ is supported in even dimensions.
Moreover, if $\Gamma$ has no bad points, then $\HFpRed(-\partial W(\Gamma))=0$.
\end{theorem}

\begin{proof}
These statements are Corollary 1.4 and
Lemma 2.6 respectively
from \cite{SomePlumbs}.
\end{proof} 

We state the following result in the case of either integral surgeries
or $1/q$ surgeries for integral $q$. However, a corresponding
statement can also be proved for $p/q$ surgeries by adapting arguments
from \cite{KMOS}, where an analogous result is proved using the
Seiberg--Witten monopole equations. We do not pursue this generalization
here, however.

\begin{theorem}
\label{thm:HFpObstruction}
Let $K\subset S^3$ be a knot in the three--sphere, and suppose that
there is an integer $p\geq 0$ and a negative definite or
semi-definite graph $\Gamma$
with only one bad point with the property that 
$$S^3_p(K)\cong -\partial W(\Gamma),$$ 
then all the elements of $\HFpRed(S^3_0(K))$ have odd
$\Zmod{2}$--grading.  Similarly, if $p>0$ and 
$$S^3_{1/p}(K)\cong -\partial W(\Gamma),$$
all the elements of a twisted Floer homology group
$\uHFpRed(S^3_0(K),\Q[\Zmod{p}])$ have odd $\Zmod{2}$--grading.  In
particular, in either case, all the torsion coefficients $t_i(K)$
(cf equation \eqref{eq:DefTorsion}) are non-negative.
\end{theorem}

\begin{proof}
If $\Gamma$ is a weighted graph, we let $Y(\Gamma)$ denote the oriented 
three--manifold obtained as $\partial W(\Gamma)$.

Assume first that $p\neq 0$, and $S^3_p(K)\cong -Y(\Gamma)$.  Consider now
the integer surgeries long exact sequence,
Theorem \ref{thm:IntegerSurgeries}
(cf \cite[Theorem 9.19]{HolDiskTwo}),
according to which for each $i\in\Zmod{p}$, we have the long exact sequence
\begin{multline*}
\ldots \stackrel{c}\longrightarrow \HFp(S^3)
\stackrel{a}{\longrightarrow}
\bigoplus_{\{j\in\Z \big|j\equiv i\!\!\!\!\pmod{p}\}}\HFp(S^3_0(K),j) \\
\stackrel{b}{\longrightarrow}
\HFp(S^3_p(K),i) \stackrel{c}{\longrightarrow}\ldots
\end{multline*}
where the maps $a$ and $c$ preserve the absolute $\Zmod{2}$--grading,
and $b$ reverses it. We assume that $S^3_p(K)\cong -Y(\Gamma)$.
Consider an element $\xi\in\HFp(S^3_0(K))$ whose
absolute $\Zmod{2}$--grading is even. Then it follows from
Theorem \ref{thm:SeifertCalc} that $b(\xi)=0$, hence that
$\xi\in\Image a$, from which it follows at once that $\xi$ comes from
$\HFinf(S^3_0(K))$, in particular, it maps trivially to
$\HFpRed(S^3_0(K))=\HFp(S^3_0(K))/\HFinf(S^3_0(K))$.
Thus, $\HFpRed(S^3_0(K))$ is supported in odd grading, as claimed.

For $1/p$--surgeries with $p\neq 0$, we can repeat the above argument,
only now using the fractional surgeries long exact sequence,
Theorem \ref{thm:ExactFrac}
(cf \cite[Theorem 9.14]{HolDiskTwo}),
which now reads
$$\ldots  \stackrel{{c}}{\longrightarrow} 
  \HFp(S^3)  \stackrel{{a}}{\longrightarrow}  \uHFp(S^3_0(K),\Q[\Zmod{p}])
  \stackrel{{b}}{\longrightarrow}  \HFp(-Y(\Gamma))
  \stackrel{{c}}{\longrightarrow} \ldots$$
For the statement of the result where $p=0$, we extend the methods of
\cite[Section 2]{SomePlumbs} to prove that for
negative semi-definite graphs $\Gamma$ with at most one bad point,
$\HFpRed(-Y(\Gamma))$ is supported in odd degrees.  To this end,
if $v$ is a vertex in a marked graph $\Gamma$, we let $\Gamma_{-1}(v)$
denote the marked graph which agrees with $\Gamma$, except 
that the weight of $v$ for $\Gamma_{-1}(v)$ is one less than the weight of $v$
for $\Gamma$. Now,  the surgery long exact
sequence, in the form it is used in
\cite[Proposition 2.8]{SomePlumbs},
gives:
$$\ldots  \stackrel{{c}}{\longrightarrow} 
  \HFp(-Y(\Gamma-v)) \stackrel{{a}}{\longrightarrow}  \HFp(-Y(\Gamma))
  \stackrel{{b}}{\longrightarrow}  \HFp(-Y(\Gamma_{-1}(v))
  \stackrel{{c}}{\longrightarrow} \ldots$$
and again, $a$ and $c$ preserve $\Zmod{2}$--grading while $b$
reverses it.  It is straightforward to see now that $\Gamma_{-1}(v)$ is a
negative-definite graph with at most one bad point, and hence
Theorem \ref{thm:SeifertCalc} applies to it. Thus, if
$\xi\in\HFpEv(S^3_0(K))$, then its image under $b$ is trivial, and
hence $\xi=a(\eta)\in\HFp(-Y(\Gamma-v))$. Indeed, since $\Gamma-v$ has no bad
points, Theorem \ref{thm:SeifertCalc} ensures that $\HFpRed(-Y(\Gamma-v))=0$,
ie $\eta$ comes from $\HFinf(-Y(\Gamma-v))$, and hence $\xi$ comes from
$\HFinf(-Y(\Gamma))$, proving the claim.

For the statement about the Alexander polynomial we appeal to Theorem
\ref{thm:Euler}, according to which when $i\neq 0$,
$$\chi(\HFp(S^3_0(K),i))=-t_i(K)$$
(in this case
$\HFp(S^3_0(K),i)=\HFpRed(S^3_0(K),i)$); while in the case where
$i=0$, we still have that $$\chi(\HFpRed(S^3_0(K),0))\geq -t_0(K).$$

For the fractional surgeries case, we use the fact that
$$\chi(\uHFp(Y_0,i,\Q[\Zmod{p}]))=-p\cm t_i(K)$$
when $i\neq0$ (this is a combination of equation \eqref{eq:TwistedChi} and
Theorem \ref{thm:Euler}), and
$$\chi(\uHFpRed(S^3_0(K),0,\Q[\Zmod{p}]))\geq -p\cm t_0(K)$$
(this follows from equation \eqref{eq:RefinedEulerTwist}, together with
inequality \eqref{eq:DIneq}).
\end{proof}

In the case of negative Seifert orientation, we have the following:

\begin{prop}
\label{prop:NegativeCase}
Let $q>0$ be an integer, and suppose that
$S^3_{1/q}(K)$ is a Seifert fibered space with negative
Seifert orientation. Then all the elements of
$\uHFpRed(S^3_0(K),\Q[\Zmod{q}])$ have even $\Zmod{2}$--grading. Moreover,
\begin{eqnarray*}
d_{1/2}(S^3_0(K))=\OneHalf
&{\text{and}}& d_{-1/2}(S^3_0(K))=-\OneHalf.
\end{eqnarray*}
\end{prop}

\begin{proof}
Consider first for notational simplicity the case where $q=1$. 

According to \cite[Proposition 2.5]{HolDiskTwo},
if $Y$ is any rational homology three--sphere,
then $\HFpRed(Y) \newline \cong \HFpRed(-Y)$, under a map which reverses the
$\Zmod{2}$--grading. Combining Theorem \ref{thm:SeifertCalc} with the
surgery long exact sequence we see that $\HFpRedEv(S^3_{1}(K))=0$. It follows
in turn from this together with the surgery long exact sequence that
$\HFp(S^3)$ injects into $\HFp(S^3_0(K))$ and hence again by the exact
sequence that the natural map from $\HFp_{\odd}(S^3_0(K))$ to
$\HFp_{\ev}(S^3_1(K))$ induces an isomorphism and in particular that
$\HFp_{\red,\odd}(S^3_0(K))=0$.

For all $i\neq 0$, $\HFp_{\red}(S^3_0(K),i)=\HFp(S^3_0(K))$. Thus, we
can conclude from Theorem~\ref{thm:Euler} that for all $i\neq 0$,
$t_i(K)<0$.

Indeed, the injectivity of $\HFp(S^3)$ in $\HFp(S^3_0(K))$ ensures
$d_{-1/2}(S^3_0(K))=-1/2$.  This implies that
\begin{equation}
\label{eq:BoundAbove}
d_{1/2}(S^3_0(K))\leq 1/2,
\end{equation}
according to inequality \eqref{eq:DIneq}. Now, since the
map from $\HFp_{\odd}(S^3_0(K),0)$ to $\HFp_{\ev}(S^3_1(K))$ drops degree
by $1/2$ we see that $d(S^3_1(K))\leq
d_{1/2}(S^3_{0}(K))-1/2$. But a negatively oriented Seifert fibered space
bounds a smooth, negative-definite four--manifold (given by its
plumbing description), so we can apply Theorem \ref{thm:NoIntForm} to
conclude that $d_{1/2}(S^3_0(K))\geq 1/2$. Putting this together with
equation \eqref{eq:BoundAbove}, we conclude that $d_{1/2}(S^3_0(K))=1/2$.

In the case where $q>1$, we use the long exact sequence
for fractional surgeries (Theorem \ref{thm:ExactFrac})
to obtain the corresponding statement for twisted coefficients. Note that
the correction terms are independent of the choice of $q$, according to
Lemma \ref{lemma:IndepTwist}.
\end{proof}

{\bf Proof of Theorem \ref{intro:AlexanderObstruction}}\qua
Assume $S^3_{1/q}(K)$ is Seifert fibered. By reflecting $K$ if necessary,
we can assume that $q>0$. 

When $S^3_{1/q}(K)$ is a positively oriented Seifert space,
Theorem \ref{thm:HFpObstruction} shows that $\uHFpRed(S^3_0(K))$ has
odd parity. We then use Theorem \ref{thm:Euler} to conclude that all
the torsion coefficients of $K$ are non-negative
(equation \eqref{eq:TwistedChi} for $i\neq 0$ and
equation \eqref{eq:RefinedEulerTwist} for $i=0$).  In the case where
$S^3_{1/q}(K)$ is a negatively oriented Seifert fibered space, we use
Proposition \ref{prop:NegativeCase}, to conclude that
$\uHFpRed(S^3_0(K),\Q[\Zmod{p}])$ is supported in even parity.  From
Proposition \ref{prop:TwistedGrows} we conclude that
$\uHFpRed(S^3_0(K))$ is supported in even parity as well. The Euler
characteristic relations (Theorem \ref{thm:Euler}, together with
equations \eqref{eq:TwistedChi} and \eqref{eq:RefinedEulerTwist})
now imply that all the torsion coefficients are non-positive in this case.
\endproof

\section{Relationship with the knot Floer homology}
\label{sec:KnotInvariants}

Theorem \ref{intro:HFKObstruction} follows from
Theorem \ref{thm:HFpObstruction}, together with the relationship
between $\HFKa(S^3,K)$ and the Heegaard Floer homology of surgeries
along $K$ developed in \cite[Section 4]{Knots},
see also \cite{RasmussenThesis}.

In \cite[Corollary 4.5]{Knots}, it is
shown that if $\HFKa(S^3,K,i)=0$ for all $i>d$, then 
$\HFKa(S^3,K,d)\cong
\HFp(S^3_0(K),d-1)$
as relatively $\Zmod{2}$--graded Abelian
groups. (With the conventions of \cite{HolDiskTwo},
the isomorphism here reverses the $\Zmod{2}$--grading.)
In the following lemma, we give the relevant statement when $d=1$.

\begin{lemma} 
\label{lemma:RelateLemma}
If $\HFKa(S^3,K,i)=0$ for all $|i|>1$, then 
$$\HFpRedEv(S^3_0(K),0)\cong \HFKa_{\odd}(S^3,K,1).$$
\end{lemma}

\proof
We use
\cite[Theorem 4.4]{Knots}.
According to that theorem, we have a $\Z\oplus \Z$--filtered complex
$C=\CFKinf(S^3,K,0)$, which admits quotient complexes
\begin{eqnarray*}
C\{i\geq 0~\text{or}~j\geq d-1\}&\sim& \CFp(S^3_p(K),[d-1]), \\
C\{i\geq 0\}&\sim& \CFp(S^3)
\end{eqnarray*}
where here $\sim$ denotes relatively graded, absolutely $\Zmod{2}$--graded
chain homotopy equivalence, $C\{i\geq 0, j\geq 0\}$ denotes the
quotient complex of $C$ by all elements whose filtration level $(i,j)$
has both $i<0$ and $j<0$, $p$ is any sufficiently large positive
integer, and $[d-1]$ is a $\SpinC$ structure over $S^3_p(K)$, which is
naturally $\SpinC$--cobordant to a $\SpinC$ structure over
$S^3_0(K)$ whose first Chern class is $2(d-1)$ times a generator
of $H^2(S^3_0(K);\Z)$.
For more on this, see \cite[Section 4]{Knots}.
Clearly, we have the short exact sequence
$$0 \longrightarrow  C\{i<0~{\text{and}}~j\geq d-1\}  \longrightarrow  
C\{i\geq 0~\text{or}~j\geq d-1\}
 \longrightarrow 
C\{i\geq 0\}
 \longrightarrow 0$$
The hypothesis that $\HFKa(S^3,K,j)=0$ for all $j>d$ ensures, by
taking filtrations, that $$H_*(C\{i<0~{\text{and}}~j\geq
d-1\})\cong\HFKa(S^3,K,d).$$ Thus, we obtain the long exact sequence
$$\ldots  \longrightarrow  \HFKa(S^3,K,d)  \longrightarrow 
  \HFp(S^3_p(K),[d-1])  \longrightarrow \HFp(S^3)
  \stackrel{{\delta}}{\longrightarrow} \ldots$$
As the above maps are $U$--equivariant, it is easy to see (from the
$\Z[U]$--module structure of $\HFp(S^3)$) that $\delta$ is the trivial
map, thus $\HFp(S^3_p(K),[d-1])$ contains a $\Z[U]$--submodule isomorphic
to $\HFKa(S^3,K,d)$. In fact, we have that
$$\HFKa_{\odd}(S^3,K,d)\cong \HFpOdd(S^3_p(K),[d-1]).$$
Specializing to the case where $d=1$, we consider the integral surgeries long
exact sequence. This gives exactness for
$$\HFp_{\ev}(S^3)  \longrightarrow \HFpEv(S^3_0(K),0)  \longrightarrow
  \HFpOdd(S^3_p(K),[0]) \longrightarrow 0.$$
It is easy to see that the image of $\HFp_{\ev}(S^3)$ inside
$\HFpEv(S^3_0(K),0)$ coincides with the image of $\HFinf_{\ev}(S^3_0(K),0)$
inside $\HFpEv(S^3_0(K),0)$. Thus, 
$$\HFpRedEv(S^3_0(K),0)\cong
  \HFpRedOdd(S^3_p(K),[0])=\HFpOdd(S^3_p(K),[0]) \eqno{\qed}$$

\noindent{\bf {Proof of Theorem \ref{intro:HFKObstruction}.}}
According to
\cite[Corollary 4.5]{Knots},
if $\HFKa(S^3,K,i)=0$ for all $i>d>1$, then
\begin{equation}
\label{eq:Relate}
\HFKa(S^3,K,d)\cong \HFp(S^3_0(K),d-1)
\end{equation}
under an isomorphism which reverses parity.

Suppose that $S^3_{1/p}(K)$ is a positively oriented Seifert fibered
space. According to Theorem \ref{thm:HFpObstruction},
$\uHFpRed(S^3_0(K);\Q[\Zmod{p}])$ is supported in odd degrees.  From
this and Proposition \ref{prop:TwistedGrows}, we can conclude the same
for $\HFpRed(S^3_0(K))$.  When the genus of $K$ is greater than one,
we have that
$$\HFp(S^3_0(K),g-1)=\HFpRed(S^3_0(K),g-1),$$
and hence all its elements have odd parity. 
Thus, in view of 
equation \eqref{eq:Relate}, $\HFKa(S^3,K,g)$ is supported in even
degrees (and it is non-trivial, cf equation \eqref{eq:DetectsGenus}). 
The case where $g=1$ follows from a similar argument, using 
Lemma \ref{lemma:RelateLemma} in place of equation \eqref{eq:Relate}

Suppose that $g>1$ and $S^3_{1/p}(K)$ is a negatively oriented Seifert
fibered space. Now, according to Proposition \ref{prop:NegativeCase},
we conclude that $\uHFpRed(S^3_0(K);\Q[\Zmod{p}])$ is supported entirely
in even degrees, and hence according to Proposition \ref{prop:TwistedGrows},
it follows that  $\HFpRed(S^3_0(K))$ is supported entirely in even degrees. 
Since $g>1$, the same can be said about $\HFp(S^3_0(K),g-1)$.
From equation \eqref{eq:Relate} it now follows that $\HFKa(K,g)$
is supported in odd degrees (and it is again non-trivial according
to equation \eqref{eq:DetectsGenus}).
\endproof

\section{Surgeries giving $\Sigma(2,3,5)$ and $\Sigma(2,3,7)$}

In general, if $Y$ is a three--manifold and $K\subset S^3$ is a knot with 
the property that $S^3_r(K)\cong Y$, then the exact sequences determine
the Floer homology of the zero--surgery only up to some indeterminacy.
However, in cases where the Floer homology of $Y$ has sufficiently 
small rank, this indeterminacy is eliminated.

{\bf Proof of Theorem \ref{thm:PoincareSphere}}\qua
From Casson's invariant, it follows that if
$S^3_{1/q}(K)=\Sigma(2,3,5)$, then $q=\pm 1$. Indeed, by reversing
orientation if necessary, we have a knot with $S^3_{\pm
1}(K)=-\Sigma(2,3,5)$. It follows from Donaldson's diagonalization
theorem that the surgery coefficient must be $+1$.  Next, we apply the
long exact sequence for $+1$ surgeries, together with the
fact that $\HFp(-\Sigma(2,3,5))\cong \InjMod{}$, with
$d(\Sigma(2,3,5))=-2$.  It follows easily that $\HFp(S^3_0(K))$ is
uniquely determined from the fact that $S^3_1(K)=-\Sigma(2,3,5)$
(compare also
\cite[Proposition 8.1]{AbsGraded}).
Indeed, we see that $\HFp(S^3_0(K))\cong
\InjMod{}\oplus\InjMod{}$, with $d_{-1/2}(S^3_0(K))=-1/2$ and
$d_{1/2}(S^3_0(K))=-3/2$. In particular, for all $i\neq 0$,
$\HFp(S^3_0(K),i)=0$. This establishes that the Floer homology of
$S^3_0(K)$ is the same as the Floer homology of $S^3_0(T_{2,3})$,
where $T_{2,3}$ denotes the right-handed trefoil, and in particular,
so is its Alexander polynomial. The fact that the knot Floer homology
coincides with that of the trefoil follows from the main result
of~\cite{NoteLens} (cf
\cite[Theorem 1.2]{NoteLens}), according
to which if $K\subset S^3$ satisfies
the property that $\HFpRed(S^3_p(K))=0$ for some integer $p$, then
$\HFKa$ is uniquely determined by the Alexander polynomial of $K$.
The remark about the Seifert genus follows from
equation \eqref{eq:DetectsGenus}
(\cite[Theorem 1.2]{HolDiskGenus}).
\endproof

{\bf Proof of Theorem \ref{thm:Bries2}}\qua
Again, from Casson's invariant it follows at once that if
$S^3_r(K)=\Sigma(2,3,7)$, then $r=\pm 1$. In the case where $S^3_{-1}(K)\cong
\Sigma(2,3,7)$, a chase of the surgery long exact sequence shows that
$\HFp(S^3_0(K))$ coincides with that for the the right-handed trefoil;
and in particular $\HFp(S^3_0(K),i)=0$ for all $i\neq 0$. Now, by
equation \eqref{eq:Relate}, it follows that $\HFKa(S^3,K,d)=0$ for all
$d>1$, and hence, by equation \eqref{eq:DetectsGenus}, the genus of
$K$ is one.

In the case where $S^3_{+1}(K)\cong \Sigma(2,3,7)$, a chase of the
surgery long exact sequence once again yields that $\HFp(S^3_0(K))$ is
uniquely determined, in particular, it is isomorphic to
$\HFp(S^3_0(K_0))$, where here $K_0$ is the figure eight knot. More
explicitly, $\HFp(S^3_0(K))=\bigoplus_{i\in\Z}\HFp(S^3_0(K),i)\cong
\InjMod{}\oplus\InjMod{}\oplus \Q$, (where the last summand has even
parity) and $d_{\pm 1/2}(S^3_0(K))=\pm 1/2$. 
We complete the argument as before.
\endproof


\Addresses\recd
\end{document}